\documentclass[12pt, abstracton]{scrartcl}

\usepackage[margin=2.5cm]{geometry}

\usepackage{pstricks, pst-tree, pst-node} 
\usepackage[latin1]{inputenc}
\usepackage{amssymb}
\usepackage{xfrac}     
\usepackage{lmodern}
\usepackage{stmaryrd}
\usepackage{graphicx}
\usepackage{amsmath}  
\usepackage{amsthm}
\usepackage{wrapfig} 
\usepackage{color}  
\usepackage{microtype} 
\usepackage{setspace}
\usepackage{amsfonts}
\usepackage{wasysym}
\usepackage{cancel}
\usepackage{verbatim}
\usepackage{dsfont}
\usepackage{paralist}

\DeclareMathOperator{\R}{\mathbb{R}}
\DeclareMathOperator{\Z}{\mathbb{Z}}
\DeclareMathOperator{\Q}{\mathbb{Q}}

\DeclareMathOperator{\N}{\mathbb{N}}

\newcommand{\qf}[1]{\langle #1\rangle}
\newcommand{\pyth}[1]{#1_\text{py}}
\newcommand{\laurent}[1]{(\!(#1)\!)}
\newcommand{\Pfister}[1]{\langle\!\langle #1\rangle\!\rangle}

\newtheoremstyle{plain2}
  {10pt}   % ABOVESPACE
  {10pt}   % BELOWSPACE
  {\itshape}  % BODYFONT
  {0pt}       % INDENT (empty value is the same as 0pt)
  {\bfseries} % HEADFONT
  {}         % HEADPUNCT
  {5pt plus 1pt minus 1pt} % HEADSPACE
  {}          % CUSTOM-HEAD-SPEC
  
 \newtheoremstyle{beweis}
  {10pt}   % ABOVESPACE
  {10pt}   % BELOWSPACE
  {\normalfont}  % BODYFONT
  {0pt}       % INDENT (empty value is the same as 0pt)
  {\bfseries} % HEADFONT
  {:}         % HEADPUNCT
  {5pt plus 1pt minus 1pt} % HEADSPACE
  {}          % CUSTOM-HEAD-SPEC

\newtheoremstyle{definition2}
  {10pt}   % ABOVESPACE
  {10pt}   % BELOWSPACE
  {\normalfont}  % BODYFONT
  {0pt}       % INDENT (empty value is the same as 0pt)
  {\bfseries} % HEADFONT
  {}         % HEADPUNCT
  {5pt plus 1pt minus 1pt} % HEADSPACE
  {}          % CUSTOM-HEAD-SPEC

\theoremstyle{plain2}
\newtheorem{Satz}{Satz}[section]
\newtheorem{Lemma}[Satz]{Lemma}
\newtheorem{Proposition}[Satz]{Proposition}
\newtheorem{Theorem}[Satz]{Theorem}
\newtheorem{Korollar}[Satz]{Corollary}
\newtheorem{Frage}[Satz]{Question}

\theoremstyle{definition2}
\newtheorem{Definition}[Satz]{Definition}
\newtheorem{Bemerkung}[Satz]{Remark}
\newtheorem{Beispiel}[Satz]{Example}

\theoremstyle{beweis}
\newtheorem*{Beweis}{Proof}

\begin{document}

\begin{titlepage}
	\vspace*{\fill}
	\centering
	{\scshape\LARGE\bfseries Supreme Pfister Forms \par}
	\vspace{1.5cm}
	{\scshape Nico Lorenz\par}
	\vspace{1.5cm}
	{Fakult\"at f\"ur Mathematik, Technische Universit\"at Dortmund, D-44221 Dortmund, Germany\par}
	\vspace{2cm}
	{e-mail: nico.lorenz@tu-dortmund.de\par}
% Bottom of the page
	\vspace{2cm}
	{\today\par}
	\vspace*{\fill}
\end{titlepage}

\begin{abstract}
	We study formally real, non-pythagorean fields which have an anisotropic torsion form that contains every anisotropic torsion form as a subform. We obtain consequences for certain invariants and the Witt ring of such fields and construct examples. We obtain a theory analogous to the theory of \textit{supreme Pfister forms} introduced by Karim Becher and see examples in which the Pythagoras number for formally real fields behaves like the level for nonreal fields.
\end{abstract}


\begin{section}{General Notation}

Throughout the text, let $F$ be a formally real non-pythagorean field, i.e. a field in which $-1$ cannot be written as a sum of squares and such that there is a sum of squares that is not a square itself. In particular, we have $\mathrm{char}(F)=0$. By a quadratic form or just form for short, we always mean a nondegenerate and finite-dimensional quadratic form.\\
For a (necessarily nonreal) field $K$ of characteristic $\neq2$, K. Becher investigated the impact of the existence of an anisotropic quadratic form containing every anisotropic form as a subform on the Witt ring $WK$ and on several field invariants of $K$. Such forms are called \textit{supreme Pfister forms}. Since some invariants such as the level or the classical $u$-invariant that are considered in Becher's article do not give meaningful information for real fields, it is natural to ask whether there is an analogous concept for formally real fields. In this article, we propose the concept of \textit{supreme torsion forms}, i.e. anisotropic torsion forms that contain every anisotropic torsion form as a subform to get a resembling theory that is similarly fruitful. We obtain comparable results for the standard generalization of the $u$-invariant to formally real fields and for the Pythagoras number instead of the level.

Before giving a more detailed exposition of the content, we will fix some further notation. We will consider two forms $\varphi_1,\varphi_2$ as equal if they are isometric, in symbols $\varphi_1\cong\varphi_2$. In particular, uniqueness of a quadratic form fulfilling some properties will implicitly mean that it is unique up to isometry. In abuse of notation, we will denote the Witt class of a quadratic form $\varphi$ again by $\varphi$. By an $n$\textit{-fold Pfister form} for some $n\in\N$, we will mean a form of the shape $\Pfister{a_1,\ldots, a_n}:=\qf{1, -a_1}\otimes\qf\ldots\otimes\qf{1,-a_n}$ with $a_1,\ldots, a_n\in F^\ast$. The set of $n$-fold Pfister forms is denoted by $P_nF$, the set of forms that are similar to some $n$-fold Pfister form is denoted by $GP_nF$. Both sets generate the $n$-th power of the \textit{fundamental ideal} $IF$, which we denote by $I^nF$, both as an additive group and as an ideal.\\
By $W_tF$, we will mean the torsion part of the Witt ring and we will use the shortcut $I^n_tF:=I^nF\cap W_tF$.\\
We will further need the $u$-invariant defined by
$$u(F):=\sup\{\dim(\varphi)\mid\varphi\text{ is an anisotropic torsion form over }F \}\in\N\cup\{\infty\}.$$
We will refer to the set of nonzero sums of squares as $\sum F^{\ast2}$. The least integer $n$ such that any element in $\sum F^{\ast2}$ is a sum of $n$ squares if such an $n$ exists or $\infty$ otherwise is called the \textit{Pythagoras number} of $F$ and it is denoted by $p(F)$.
All further notation will be consistent with Lam's book \cite{Lam2005}.

In Section 2 we will formally introduce the notation of the most important forms of this article, the \textit{supreme torsion forms}. We will give easy examples of fields that have a supreme torsion forms and collect some basic properties of such fields and the supreme torsion forms itself. We will see that supreme torsion forms are always Pfister forms and that they are unique. 

Section 3 will deal with the yet to be defined concept of \textit{2-real-maximality}. We will need some preliminary results to see the fruitfulness of this concept that further allows us to obtain a necessary condition for a form to be the supreme torsion form.
We will use 2-real-maximality frequently in Section 4 in order to construct fields with a supreme torsion form that have further additional properties. 

In the final Section of this article, we will return to the question how the existence of a supreme torsion forms influences some field invariants that concern quadratic form theory. We will focus on the cardinality of $F^\ast/F^{\ast2}$ and the Pythagoras number $p(F)$.

\end{section}

\begin{section}{Introduction to Supreme Torsion Forms}

We will now give the definition of what we call \textit{supreme torsion forms} in this text. In abuse of notation we call a quadratic form over $F$ torsion if its Witt class is torsion in $WF$. 

\begin{Definition}
	Let $\varphi$ be a quadratic form over $F$. We call $\varphi$ a \textit{supreme torsion form} if $\varphi$ is an anisotropic torsion form such that every anisotropic torsion form over $F$ is similiar to a subform of $\varphi$.
\end{Definition}

\begin{Beispiel}\label{BspSupTorsForm}
	\begin{enumerate}[(a)]
		\item\label{BspSupTorsFormA} We consider a formally real field $F$ with square class group $\{1,-1, 2, -2\}$. Such a field exists by \cite[II. Remark 5.3]{Lam2005} and can be constructed by a modification of the Gross-Fischer construction. For the Witt group we have $WF\cong\Z\oplus\Z/2\Z$ and the only nontrivial torsion form is given by $\Pfister2=\qf{1,-2}$ which then trivially is a supreme torsion form.
		\item Confer \cite[Chapter II, Example 5.4 and page 45]{Lam2005}: Let $F=\Q_3\cap R$, where $R$ is a real closed subfield of the algebraic closure $\overline{\Q_3}$ of the $3$-adic numbers. Then $F$ has square class basis $\{-1, 2, 3\}$ and $3$ is not a sum of two squares. The torsion subgroup $W_tF$ is given by
		$$
		\{\mathbb H, \qf{1, -2}, \qf{1, -3},\qf{2, -6},\qf{1, -6}, \qf{-3, 6},\qf{2,-3}, \qf{1, 1, -3, -3}\}.$$ 
		One readily checks that $\qf{1, 1, -3, -3}$ is a supreme torsion form.
		\item Confer \cite[Chapter II, Example 5.7, and page 46]{Lam2005}: Analogously to the above example we can take $F=\Q_5\cap R$ where $R$ is a real closed subfield of the algebraic closure $\overline{\Q_5}$ of the $5$-adic numbers. Then $F$ has square class basis $\{-1, 2, 5\}$. The torsion subgroup $W_tF$ is given by
		 \begin{align*}
			\{\mathbb H, \qf{1, -2}, \qf{1, -5},\qf{1, -10},\qf{2, -5}, \qf{2, -10}, \qf{5, -10}, \qf{1, -2, 5, -10}\}
		\end{align*}
		and $\qf{1, -2, 5, -10}$ is a supreme torsion form as can be readily verified.
	\end{enumerate}
\end{Beispiel}

\begin{Lemma}\label{STFuInv}
	If $F$ admits a supreme torsion form $\varphi$, then $F$ has finite $u$-invariant $u(F)=\dim\varphi$.
\end{Lemma}
\begin{Beweis}
	As every anisotropic torsion form over $F$ is similar to a subform of $\varphi$, we have $u(F)\leq\dim\varphi$. But $\varphi$ is anisotropic by definition of a supreme torsion form, so we have $u(F)=\dim\varphi$. 
\end{Beweis}

\begin{Proposition}\label{STFisPfister}
	Let $\varphi$ be a supreme torsion form. Then $\varphi$ is a Pfister form and every anisotropic torsion form is a subform of $\varphi$. In particular $\varphi$ is unique up to isometry.
\end{Proposition}
\begin{Beweis}
	We choose an anisotropic torsion Pfister form $\pi$ of maximal dimension. Such a form exists as $F$ is not pythagorean so there are two dimensional torsion forms, and every torsion form of dimension greater than $\dim\varphi$ has to be isotropic according to \ref{STFuInv}. Then for every $a\in F^\ast$, $\pi\otimes\Pfister{a}$ is a torsion form. Because of the maximality of $\pi$, it is isotropic and therefore as a Pfister form hyperbolic. That means we have $a\pi\cong\pi$, i.e. $a\in G_F(\pi)=D_F(\pi)$. Thus $\pi$ is universal so that $\pi$ cannot be similar to a proper subform of any anisotropic form. Thus $\pi$ is similar to $\varphi$. Because $\pi$ is universal and round as a Pfister form, it is therefore isometric to $\varphi$. So if $\psi$ is an anisotropic torsion form and $a\in F^\ast$ such that $a\psi\subseteq\varphi$, we have $$\psi\cong a^2\psi\subseteq a\varphi\cong\varphi$$
	as $\varphi\cong\pi$ is round and universal as shown above.
\end{Beweis}

In view of the above result, we will often call a supreme torsion form \textit{the} supreme torsion form in the sequel.

\begin{Proposition}\label{STFisOnlyUniTorForm}
	Let $\varphi$ be the supreme torsion form over $F$. Then $\varphi$ is the only anisotropic universal torsion form.
\end{Proposition}
\begin{Beweis}
	Let $\psi$ be an anisotropic universal torsion form. Then $\psi$ is isometric to a subform of $\varphi$ by \ref{STFisPfister}. But as a universal form cannot be a proper subform of any anisotropic form we get $\psi\cong\varphi$.
\end{Beweis}

Considering the above proposition and \cite[3.3 Corollary]{Becher2004}, it is natural to ask the following:

\begin{Frage}\label{EinzigeUniverselleTorsionsFormIstSupremeFragezeichen}
	Let $u(F)<\infty$ and $\varphi$ be the only anisotropic universal torsion form. Is $\varphi$ the supreme torsion form?
\end{Frage}

To construct examples of every possible size we study the behaviour of supreme forms under Laurent series extensions.

\begin{Proposition}\label{STFLaurent}
	Let $\varphi\in P_nF$ be the supreme torsion form over $F$. Then $\varphi\otimes\Pfister t$ is the supreme torsion form over $F\laurent t$.\\
	If conversely $\psi=\psi_1\perp -t\psi_2$ with residue class forms $\psi_1,-\psi_2$ is the supreme torsion form over $F\laurent t$ we have $\psi_1\cong\psi_2$ and this form is the supreme torsion form over $F$.
\end{Proposition}
\begin{Beweis}
	It is clear that a form over $F\laurent t$ is torsion if and only if both residue class forms of the given form are torsion over $F$. In particular $\varphi\otimes\Pfister t$ is torsion. It is further anisotropic by a consequence of Springer's Theorem, see \cite[VI. Proposition 1.9]{Lam2005}. If $\psi=\psi_1\perp -t\psi_2$ is a torsion form over $F\laurent t$ we have that both $\psi_1$ and $\psi_2$ are subforms of $\varphi$ as $\varphi$ is the supreme torsion form over $F$. It is then clear that $\psi$ is a subform of $\varphi\otimes\Pfister t=\varphi\perp-t\varphi$. 
	
	For the converse let $\tau$ be an anisotropic torsion form over $F$. Then $\tau\otimes\Pfister t$ is an anisotropic torsion form over $F\laurent t$ which is therefore a subform of the supreme torsion form $\psi$. This is equivalent to $$i_W(\psi\perp-\tau\otimes\Pfister t)\geq\dim(\tau\otimes\Pfister t)=2\dim\tau.$$ On the other hand we have
	\begin{align*}
		i_W(\psi\perp-\tau\otimes\Pfister t)&=i_W(\psi_1\perp-\tau\perp -t(\psi_2\perp\tau))\\
		&=\underbrace{i_W(\psi_1\perp-\tau)}_{\leq\dim\tau}+\underbrace{i_W(\psi_2\perp\tau)}_{\leq\dim\tau}\leq 2\dim\tau,
	\end{align*}
	which finally yields
	$$i_W(\psi_1\perp-\tau)=\dim\tau=i_W(\psi_2\perp\tau).$$
	This implies $\tau\subseteq\psi_1$ and $-\tau\subseteq\psi_2$. As $\tau$ is an arbitrary anisotropic torsion form we see that both $\psi_1$ and $\psi_2$ are supreme torsion forms over $F$. As supreme torsion forms are unique up to isometry by \ref{STFisPfister}, we finally see $\psi_1\cong\psi_2$.
\end{Beweis}

\begin{Beispiel}
	Combining \ref{BspSupTorsForm} with \ref{STFLaurent} we can construct a field $F_n$ with an $n$-fold Pfister form as a supreme torsion form for any $n\in\N$. To do so let $F_1$ be the field as in \ref{BspSupTorsForm} \ref{BspSupTorsFormA} and define $F_n$ to be $F_n:=F_1\laurent{t_1}\cdots\laurent{t_{n-1}}$ for $n\geq 2$.
\end{Beispiel}

\begin{Proposition}
	Let $n\in\N$ be an integer and $\pi$ the supreme torsion $n$-fold Pfister form over $F$. We then have:
	\begin{enumerate}[(a)]
		\item\label{DontYouActFunny} $I_t^{n+1}F$ is trivial;
		\item $\pi$ is the unique anisotropic torsion $n$-fold Pfister form over $F$;
		\item Every anisotropic form in $I_t^{n-1}F$ is either isometric to $\pi$ or similar to an ${(n-1)}$-fold Pfister form;
		\item If $\varphi$ is a nonhyperbolic torsion form, there is a Pfister form $\psi$ over $F$ such that $\varphi\otimes\psi$ is Witt equivalent to $\pi$. If we have moreover $\varphi\in I^kF$ and $\dim\varphi<2^{k+1}$ for some $k\in\{1,\ldots, n\}$, then $\psi$ is an $(n-k)$-fold Pfister form.
	\end{enumerate}
\end{Proposition}
\begin{Beweis}
	\begin{enumerate}[(a)]
		\item Is clear with \ref{STFuInv} and the Arason-Pfister Hauptsatz.
		\item Every anisotropic torsion $n$-fold Pfister form over $F$ is isometric to a subform of $\pi$ because of \ref{STFisPfister}. Considering the dimensions we get that such a form is isometric to $\pi$.
		\item Let $\varphi\in I^{n-1}_tF\setminus\{0\}$ be an anisotropic torsion form with $\varphi\neq\pi$. According to \ref{STFisPfister}, $\varphi$ is thus a subform of $\pi$. We can therefore write $\pi\cong\varphi\perp\psi$. As $\pi,\varphi\in I_t^{n-1}F$ we further have $\psi\in I_t^{n-1}F$. The Arason-Pfister Hauptsatz yields $\dim\varphi=\dim\psi=2^{n-1}$ and therefore $\varphi$ and $\psi$ have to be similar to $(n-1)$-fold Pfister forms respectively.
		\item As $\varphi$ is not hyperbolic, every multiple of $\varphi$ is torsion and $u(F)$ is finite, we can choose an $m$-fold Pfister form $\psi$ such that $\varphi\otimes\psi$ is not hyperbolic with $m\in\N_0$ maximal. Let $\rho:=(\varphi\otimes\psi)_{\text{an}}$ denote the anisotropic part of $\varphi\otimes\psi$. Because of the maximality of $m$, for every $a\in F^\ast$, the form $\varphi\otimes\psi\otimes\Pfister a$ is hyperbolic, which is equivalent to $\rho\cong a\rho$. This implies $\rho$ to be universal. \ref{STFisOnlyUniTorForm} implies $\rho\cong\pi$ which implies $\varphi\otimes\psi=\pi\in WF$.\\
		If we have $\varphi\in I_t^kF$ with $\dim\varphi<2^{k+1}$, we get
		$$2^n=\dim(\pi)\leq\dim(\varphi\otimes\psi)< 2^{m+k+1}$$
		which implies $m\geq n-k$. We also have $m\leq n-k$ as otherwise we would get $\varphi\otimes\psi\in I_t^{n+1}$ contradicting (\ref{DontYouActFunny}).
	\end{enumerate}
\end{Beweis}

The next result is an analogue to Kneser's Lemma, see \cite[Chapter XI., 6.5]{Lam2005}.

\begin{Proposition}\label{KnesersLemmaTorsion}
	Let $\varphi$ be a non universal torsion form and $\beta$ binary torsion form over $F$. We then have $D_F(\varphi)\subsetneq D_F(\varphi\perp\beta)$.
\end{Proposition}
\begin{Beweis}
	We have $\beta\cong\qf{a,-xa}$ for some $a\in F^\ast$ and $x\in\sum F^{\ast2}$. Thus $x=e_1^2+\ldots+e_n^2$ for some $e_1,\ldots, e_n\in F^\ast$. Now assume $D_F(\varphi)= D_F(\varphi\perp\beta)$. We first show per induction that $a(e_1^2+\ldots+e_k^2)\in D_F(\varphi)$ for all $k\in\{1,\ldots, n\}$.
	For $k=1$ this follows from the equality $D_F(\varphi)= D_F(\varphi\perp\beta)$ as clearly $a\in D_F(\beta)\subseteq D_F(\varphi\perp\beta)$. If $a(e_1^2+\ldots+e_{k-1}^2)\in D_F(\varphi)$ we have
	$$a(e_1^2+\ldots+e_k^2)=a(e_1^2+\ldots+e_{k-1}^2)+ae_k^2\in D_F(\varphi)+D_F(\beta)\subseteq D_F(\varphi\perp\beta)=D_F(\varphi)$$
	as claimed.\\
	In particular, it follows $xa=(e_1^2+\ldots+e_n^2)a\in D_F(\varphi)$ so that $D_F(\varphi\perp\beta)$ is isotropic, hence universal. As we have $D_F(\varphi)=D_F(\varphi\perp\beta)=F^{\ast2}$ it follows that $\varphi$ is universal, contradicting the hypothesis.
\end{Beweis}

\end{section} 

\begin{section}{2-real maximality}

We develop a theory similar to that in \cite[Section 3]{Becher2006} in order to obtain a necessary condition for a form to be a supreme torsion form. We would like to draw attention to the fact that the basic definitions of what is called a totally positive field extension differ, but the results will sometimes read the same. We will compare both properties in detail later in \ref{VergleichTotallyPositive}.

\begin{Definition}
	A field extension $K/F$ is called \textit{totally positive} if every ordering of $F$ extends to an ordering of $K$.
\end{Definition}

In contrast to our definition, K. Becher defines a field extension $K/F$ as totally positive if every semi-ordering of $F$ extends to a semi-ordering of $K$. As any ordering is a semi-ordering, we deal with a weaker property than K. Becher in his article.

\begin{Beispiel}\label{BspTotallyPositive}
	\begin{enumerate}[(a)]
		\item Finite extensions of odd degree are totally positive due to \cite[VIII. Corollary 7.10 (1)]{Lam2005}.
		\item It is well known that any ordering of $F$ has exactly two extensions to $F\laurent t$, one in which $t$ is positive and one in which $t$ is negative. Thus $F\laurent t/F$ is totally positive for every real field $F$. Further purely transcendental extensions are totally positive.
		\item As the function field extension $F(\varphi)$ extends an ordering $P\in X_F$ if and only if $\varphi$ is indefinite at $P$, the extension $F(\varphi)/F$ is totally positive if and only if $\varphi$ is totally indefinite, see \cite[XIII. Theorem 3.1]{Lam2005}. In particular, for $\varphi=\Pfister a$ for some $a\in F^\ast\setminus F^{\ast2}$, we obtain that $F(\Pfister a)=F(\sqrt a)$ extends a given ordering $P\in X_F$ if and only if $a\in P$.
	\end{enumerate}
\end{Beispiel}

As another easy example we remark the following for later reference. This result also justifies the used terminology.

\begin{Lemma}\label{LemTotallyPositiveQuadExt}
	Let $a\in F^\ast\setminus F^{\ast2}$. The quadratic extension $F(\sqrt a)/F$ is totally positive if and only $a$ is totally positive, i.e. we have $a\in\sum F^{\ast2}$.
\end{Lemma}
\begin{Beweis}
	By \cite[VIII. Corollary 7.10 (2)]{Lam2005} we see that $F(\sqrt a)/F$ is totally positive if and only if $a\in P$ for all $P\in X_F$. By Artins theorem, see \cite[VII. 1.12]{Lam2005}, this is equivalent to $a\in\sum F^{\ast2}$.
\end{Beweis}

As a first step we will study how the term totally positive behaves in the case of towers of field extensions.

\begin{Lemma}\label{TotallyPositiveFieldTower}
	Let $K/E/F$ be a tower of field extensions. We then have:
	\begin{enumerate}[(a)]
		\item\label{TotallyPositiveFieldTowerA} If $K/E$ and $E/F$ both are totally positive then so is $K/F$.
		\item\label{TotallyPositiveFieldTowerB} If $K/F$ is totally positive, so is $E/F$.
	\end{enumerate}
\end{Lemma}
\begin{Beweis}
	\begin{enumerate}[(a)]
		\item This is obviously true.
		\item Every ordering of $F$ can be extended to an ordering of $K\supseteq E$ by hypothesis which then can be restricted to an ordering on $E$, which clearly extends the given ordering.
	\end{enumerate}
\end{Beweis}

In order to generalize \ref{LemTotallyPositiveQuadExt} to arbitrary subfields of the pythagorean closure, we need the following lemma. It is a generalization of \cite[Proposition 1.1.10]{ZalesskiiRibes2000}.

\begin{Proposition}\label{ExtendOrderingsDirectLimit}
	Let $(F_i, f_{ij})$ be a directed system indexed by an index set $I$ of formally real fields $F_i$ with respective space of orderings $X_i\neq\emptyset$. Let $j\in I$ and $\emptyset\neq X\subseteq X_j$ such that for all $k\in I$ with $j\leq k$, every ordering in $X$ has an extension to $F_k$. Then every ordering in $X$ has an extension to $F:=\varinjlim F_i$.
\end{Proposition}
\begin{Beweis}
	It is well known that $X_F$ is isomorphic to the inverse limit $\varprojlim X_i$ where the homomorphisms in this inverse system are given by the restrictions of the orderings, denoted by $\varphi_{ij}:X_j\to X_i$ for $i,j\in I$ with $i\leq j$. \\
	Let now $I':=\{k\in I\mid k\geq j\}$ be the cofinal subset of $I$ of indices greater than or equal to $j$ (in the given ordering of $I$).
	
	For every $\alpha\in X$ and every $k\in I'$, the set $Y_k=\varphi^{-1}_{jk}(\alpha)$ is not empty by hypothesis. As the $\varphi_{jk}$ are continuous by \cite[Corollary, page 272]{Lam2005} and $\{\alpha\}\subseteq X_j$ is compact as a singleton set in a boolean space, the $Y_k$ are all nonempty and compact. As subspaces of boolean spaces, the $Y_k$ are further totally disconnected and Hausdorff. Since the $\varphi_{ik}$ are just given by restriction we clearly have $\varphi_{ik}(Y_k)\subseteq Y_i$ for all $i,k\in I'$ with $i\leq k$. We thus have an inverse subsystem $(Y_i, \varphi_{ik}, I')$ of $(X_i, \varphi_{ik}, I')$ consisting of compact Hausdorff totally disconnected topological spaces, whose inverse limit $\displaystyle \varinjlim_{I'} Y_{i'}\subseteq \varinjlim_{I'} X_{i'}$ is not empty according to \cite[Proposition 1.1.3]{ZalesskiiRibes2000}.\\
	As $I'$ is cofinal in $I$, we have a canonical isomorphism 
	$$\varinjlim_IX_i\to \varinjlim_{I'} X_{i'}$$
	by \cite[Lemma 1.1.9]{ZalesskiiRibes2000}. By the choice of the $Y_k$ every preimage of any element of $\displaystyle \varinjlim_{I'} Y_{i'}$ corresponds to an ordering of $F$ extending $\alpha$.
\end{Beweis}

\begin{Proposition}\label{TotallyPositivePyth}
	If $K\subseteq F_{\text{py}}$ is a subfield of the pythagorean closure of $F$ then $K/F$ is totally positive.
\end{Proposition}
\begin{Beweis}
	By \ref{TotallyPositiveFieldTower} (\ref{TotallyPositiveFieldTowerB}) it is enough to show that $\pyth F/F$ is totally positive. This is shown for example in \cite[Chapter VIII. Corollary 4.6]{Lam2005}. We further give an alternative proof using techniques that will reoccur in the next section.\\
	By \cite[Lemma 31.16]{ElmanKarpenkoMerkurjev2008} $\pyth F$ is the direct limit of all field extensions $E/F$ such that there is a tower of field extensions $F=F_0\subsetneq F_1\subsetneq\ldots\subsetneq F_n=E$ with $F_k=F_{k-1}(\sqrt {z_{k-1}})$, where $z_{k-1}=1+x_{k-1}^2$ for some $x_{k-1}\in F_{k-1}$. Since each $F_k/F_{k-1}$ is totally positive by \ref{LemTotallyPositiveQuadExt}, so is $E/F$ according to \ref{TotallyPositiveFieldTower} (\ref{TotallyPositiveFieldTowerA}). Since the property of being totally positive is preserved under direct limits by \ref{ExtendOrderingsDirectLimit}, the conclusion follows.
\end{Beweis}

\begin{Bemerkung}\label{VergleichTotallyPositive}
	We will now compare our version of totally positive field extensions with the strong version of totally positive field extensions introduced by K. Becher, to which we will refer as \textit{strongly totally positive}. As mentioned above, since any ordering is a semi-ordering, strongly totally positive field extensions are always totally positive. A quadratic extension is strongly totally positive if and only if it is totally positive due to \ref{LemTotallyPositiveQuadExt} and \cite[Proposition 3.2]{Becher2006}.\\
	Further, for any subfield $E$ of the pythagorean closure of $F$, the extension $E/F$ is strongly totally positive due to \ref{TotallyPositivePyth} and \cite[Corollary 3.3]{Becher2006}.
	We finally consider the case of the function field extension of an anisotropic form $\varphi$ of dimension at least 3. The extension $F(\varphi)/F$ is strongly totally positive if and only if $n\times\varphi$ is isotropic for some $n\in\N$ by \cite[Theorem 3.4]{Becher2006}. Over so-called \textit{SAP fields} any form is totally indefinite if and only if some suitable multiple is isotropic by \cite[(9.1) Theorem]{Prestel}, where the if-part obviously holds over any field. We therefore obtain that both concepts coincide over SAP fields. But by \cite[(2.12) Theorem]{Prestel}, for any SAP field $F$, there is a form of the shape $\qf{1,a, b, -ab}$ for suitable $a,b\in F^\ast$ that is totally indefinite but $n\times\qf{1,a,b,-ab}$ is anisotropic for all $n\in\N$. As a concrete example, we can consider the field $F=\R\laurent s\laurent t$ and $\varphi=\qf{1,s, t, -st}$. Then $\varphi$ is of the desired shape but clearly, for any $n\in\N$, $n\times\varphi$ is anisotropic. Thus $F(\varphi)/F$ is a totally positive field extension that is not strongly totally positive.
\end{Bemerkung}

With the following example we will disprove the reverse implications in the above results.

\begin{Beispiel}
	Let $F=\Q, E=\Q(\sqrt2), K=E(\sqrt{\sqrt2})=\Q(\sqrt[4]2)$. As $K$ is real and $F=\Q$ has a unique ordering $K/F$ is totally positive. Now \ref{TotallyPositiveFieldTower} (\ref{TotallyPositiveFieldTowerB}) implies $E/F$ to be totally positive. But $K/E$ is not totally positive as the ordering of $E$ in which $\sqrt2<0$ cannot be extended to $K$. So $K/F$ (and $E/F$) totally positive does not in general imply $K/E$ to be totally positive for a tower of field extensions $K/E/F$, the other implication in \ref{TotallyPositiveFieldTower} (\ref{TotallyPositiveFieldTowerA}) is therefore false.\\
	We further have $F_\text{py}=E_\text{py}$ and $\sqrt2$ is not totally positive in $E$. We thus see ${K\not\subseteq \pyth E=\pyth F}$ so that the other implication in \ref{TotallyPositivePyth} is false as well.
\end{Beispiel}

As we have just seen that totally positive field extensions do not behave well in general when enlarging the base field, it seems natural to refine the term. We will see that this refinement is exactly the modification we need.

\begin{Definition}
	A field extension $K/F$ is called \textit{hereditarily totally positive} if $K/E$ is totally positive for all intermediate fields $F\subseteq E\subseteq K$.
\end{Definition}

\begin{Beispiel}
		Let $E/F$ be a field extension of odd degree and $K$ an intermediate field. As we have $[E:F]=[E:K]\cdot[K:F]$, the extension $E/K$ is also of odd degree and therefore totally real, see \cite[Chapter 3, 1.10 Theorem (ii)]{Scharlau}. Therefore odd degree extensions are hereditarily totally real.
\end{Beispiel}

\begin{Proposition}\label{QuadClosHerTotallyPositive}
	Let $K\subseteq F_2$ be a subfield of the quadratic closure of $F$. Then $K/F$ is hereditarily totally positive if and only if we have $K\subseteq \pyth F$. 
\end{Proposition}
\begin{Beweis}
	We assume first $K\subseteq\pyth F$ and let $E$ be an intermediate field $K/E/F$. We then have $E\subseteq K\subseteq\pyth F=\pyth K$ and the implication then follows by \ref{TotallyPositivePyth}.\\
	For the converse we consider the field $E:=K\cap\pyth F$. Since we are done if we have $E=K$ we now assume $E\neq K$.
	As $K/E$ is an algebraic extension within the quadratic closure of $F$, there is some quadratic extension of $E$ within $K$ by \cite[Section 18, Problem 1]{Morandi}. As we have $\mathrm{char} E\neq 2$, using a well known fact in elementary field theory yields the existence of some $a\in E$ such that we have $E\subsetneq E(\sqrt a)\subseteq K$.
	By the choice of $a$ we have $\sqrt a\notin \pyth F=\pyth E$. Thus the binary form $\qf{1, -a}$ is not isometric to the hyperbolic plane over $\pyth E$, which means that $\qf{1, -a}\notin W(\pyth F/F)=W_tF$ by \cite[VIII. Theorem 4.10]{Lam2005}. By Pfister's local global principle this means that $a$ is not totally positive in $E$. So there is an ordering $P$ on $E$ in which $a$ is negative. By \ref{BspTotallyPositive} (c) $P$ does not have any extension to $E(\sqrt a)$ and therefore cannot have any extension to $K\supseteq E(\sqrt a)$. Thus $K/E$ is not totally positive, a contradiction.
\end{Beweis}

\begin{Lemma}\label{TotallyRealAlgebraic}
	If $E/F$ is hereditarily totally positive field extension then $E/F$ is algebraic.
\end{Lemma}
\begin{Beweis}
	Let $E/F$ be a field extension that is not algebraic. Then there is some $x\in E$ that is transcendental over $F$. Consider the purely transcendental field extension $K:=F(x^2)$ of $F$ which is a proper subfield of $E$. By \cite[Chapter VIII. Examples 1.13 (C)]{Lam2005} there is some ordering on $K$ in which $x^2$ is negative. Such an ordering cannot be extended to an ordering on $E$ since $x^2$ is a square in $E$ and therefore totally positive in $E$.
\end{Beweis}

\begin{Korollar}\label{HereditarilyTotallyRealFunctionField}
	Let $\varphi$ be a quadratic form over $F$ such that the function field $F(\varphi)$ is defined, i.e. $\dim\varphi\geq2$ and if $\dim\varphi=2$ we assume $\varphi$ to be anisotropic. The function field extension $F(\varphi)/F$ is hereditarily totally positive if and only of $\varphi$ is a non hyperbolic binary torsion form.
\end{Korollar}
\begin{Beweis}
	We must have $\dim\varphi=2$, $\varphi$ anisotropic since otherwise, $F(\varphi)/F$ is not an algebraic extension contradicting \ref{TotallyRealAlgebraic}. We therefore have $\varphi\cong b\qf{1,-a}$ for some ${a,b\in F^\ast}$ and $F(\varphi)=F(\sqrt a)$. In this case, the only proper subfield of $F(\varphi)$ containing $F$ is $F$ itself. Thus $F(\varphi)/F$ is hereditarily totally positive if and only if $F(\varphi)/F$ is totally positive. By \ref{LemTotallyPositiveQuadExt} this is the case if and only if $a$ is totally positive, which is equivalent to $\varphi$ being torsion.
\end{Beweis}

As a last step in this section, we would like to deduce a necessary condition for a quadratic form to be a supreme torsion form using the above theory of hereditarily totally positive field extensions.

\begin{Definition}
	Let $\varphi$ be a form over $F$ and $E/F$ a formally real non-pythagorean field extension. Then $\varphi$ is called \textit{2-real-maximal} over $E$ if $\varphi_E$ is anisotropic but $\varphi_K$ is isotropic for every nontrivial hereditarily totally positive $2$-extension $K/E$. 
\end{Definition}

\begin{Proposition}\label{Construct2RealMaximal}
	Let $\varphi$ be an anisotropic torsion form over $F$. Then there exists a hereditarily totally positive 2-extension $K/F$ such that $\varphi$ is 2-real-maximal over $K$.
\end{Proposition}
\begin{Beweis}
	We define $\mathcal C$ to be the class $$\mathcal C:=\{K\subseteq \pyth F\mid K/F\text{ is a field extension, }\varphi_K\text{ is anisotropic}\}.$$ We clearly have $F\in\mathcal C$. Further every chain in $\mathcal C$ has its union as an upper bound. Thus Zorn's lemma implies the existence of a maximal element $K\in\mathcal C$. As $\varphi$ is torsion and we require $\varphi_K$ to be anisotropic, we have $K\subsetneq \pyth F$. If now $E/K$ is a nontrivial hereditarily totally positive 2-extension we have $E\subseteq \pyth K=\pyth F$ by \ref{QuadClosHerTotallyPositive}. Since we have chosen $K\in\mathcal C$ to be maximal, $\varphi_E$ has to be isotropic, i.e. $\varphi$ is 2-real maximal over $K$.
\end{Beweis}

\begin{Bemerkung}
	In \ref{Construct2RealMaximal} we required $\varphi$ to be a torsion form to guarantee that the constructed field $K$ is not pythagorean, in line with our general assumption that we want to work only over nonpythagorean fields. Over a pythagorean field the theory of 2-real-maximality would not be fruitful as any form would trivially be 2-real-maximal. As an aside we further remark that \cite[Example 3.6 and Lemma 3.1]{Becher2006} implies the existence of a totally indefinite form over a field $F$ that is still anisotropic over $\pyth F$ (and hence no torsion form over $F$), such that this form would lead to $K=\pyth F$ in \ref{Construct2RealMaximal} when imitating the proof.
\end{Bemerkung}

\begin{Proposition}\label{Connemara}
	Let $\varphi$ be an anisotropic quadratic form over $F$. Then the following are equivalent:
	\begin{enumerate}[(a)]
		\item The form $\varphi$ is 2-real maximal over $F$;
		\item\label{PropEqu2realMax} Every anisotropic binary torsion form is similar to a subform of $\varphi$.
	\end{enumerate}
\end{Proposition}
\begin{Beweis}
	(i) $\Rightarrow$ (ii): Let $\beta$ be an anisotropic binary torsion form. Then $\beta$ is similar to $\qf{1, -a}$ for some $a\in\sum F^{\ast2}\setminus F^{\ast2}$. We therefore have $F(\sqrt a)\subseteq \pyth F$. As $\varphi$ is 2-real maximal over $F$, $\varphi_{F(\sqrt a)}$ has to be isotropic which is equivalent to $\Pfister a=\qf{1,-a}\subseteq \varphi$.\\
	(ii) $\Rightarrow$ (i): Let $K/F$ be a nontrivial hereditarily totally positive 2-extension. According to \ref{QuadClosHerTotallyPositive} we have $K\subseteq \pyth F$. As in the proof of this proposition $K$ contains a quadratic extension of $F$. That means there is an $a\in F^\ast$ such that $\qf{1,-a}$ is not isotropic over $F$ but over $K$ and therefore also over $\pyth F\supseteq K$. This means that $\qf{1,-a}\in W(\pyth F/F)=W_tF$ is a binary torsion form over $F$ and thus a subform of $\varphi$. So $\varphi_K$ is isotropic.
\end{Beweis}

\begin{Korollar}\label{SupremeImplies2RealMaximal}
	If $\varphi$ is the supreme torsion form over $F$ then $\varphi$ is 2-real-maximal over $F$.
\end{Korollar}
\begin{Beweis}
	This is clear since $\varphi$ obviously satisfies \ref{Connemara} \ref{PropEqu2realMax}.
\end{Beweis}

The other implication in \ref{SupremeImplies2RealMaximal} is not true. We will give two different reasons to disprove this: \\
Firstly, we will show that we can construct a field such that there is more than one Pfister form that is 2-real maximal over this field simultaneously, see \ref{PropConstruction}. Since supreme torsion forms are unique up to isometry by \ref{STFisPfister}, this will finish the first argument.\\
Secondly, we can give an example of a field with a 2-real maximal form that is not even a torsion form. This will be done in \ref{2RealMaximalNotTorsion}.\\

As supreme torsion forms are always Pfister forms and 2-real maximal, it is convenient to have a criterion for 2-real-maximality specifically for Pfister forms. This is done in the following using the set of represented elements of the pure part of the given Pfister form.

\begin{Proposition}
	Let $\pi$ be an anisotropic Pfister form. Then $\varphi$ is 2-real-maximal over $F$ if and only if we have $-\sum F^{\ast2}\setminus -F^{\ast2}\subseteq D_F(\pi')$.
\end{Proposition}
\begin{Beweis}
	For the if-part, let $\qf{a,b}$ be a binary anisotropic torsion form. We then have 
	$$ab\in -\sum F^{\ast2}\setminus-F^{\ast2}\subseteq D_F(\pi').$$
	We thus have $\qf{1,ab}\subseteq\pi$. As $\qf{1,ab}$ is similar to $\qf{a,b}$, this implication follows in view of \ref{Connemara}.\\
	For the opposite direction, let $a\in -\sum F^{\ast2}\setminus-F^{\ast2}$. As the form $\qf{1,a}$ is torsion and anisotropic by the choice of $a$, it is similar to a subform of $\pi$ by \ref{Connemara}, i.e. there is some $x\in F^{\ast}$ with $x\qf{1,a}\subseteq\pi$. This in particular implies ${x\in D_F(\pi)=G_F(\pi)}$, so that we even get $\qf{1,a}\subseteq\pi$. Witt Cancellation now yields $\qf{a}\subseteq\pi'$, i.e. ${a\in D_F(\pi')}$.
\end{Beweis}

\end{section}

\begin{section}{Construction of Examples}
The aim of this section is to introduce methods to construct fields with supreme torsions forms or fields that have a given set of forms that are 2-real maximal and which fulfil several additional properties in order to fill the gaps of the former section. 

\begin{Definition}
	Let $\mathcal C$ be a class of field extensions of a field $F$. The class $\mathcal C$ is called \textit{admissible} if the following holds:
	\begin{enumerate}
		\item[(AD1)]\label{admissible1} The class $\mathcal C$ is not empty.
		\item[(AD2)]\label{admissible2} The class $\mathcal C$ is closed under direct limits (in the category of field extensions of $F$).
		\item[(AD3)]\label{admissible3} If $K\in\mathcal C$ and $E\subseteq K$ is a subfield of $K$ with $F\subseteq E$ then we have $E\in\mathcal C$.
	\end{enumerate}
\end{Definition}

\begin{Bemerkung}\label{DirectSystemSubfields}
	Recall that, for a direct system $(E_i)_{i\in I}$ of field extensions of $F$ with direct limit $E=\varinjlim_{i\in I}E_i$, we have embeddings $\varphi_i:E_i\to E$. Replacing $E_i$ with $\varphi_i(E_i)$ for all $i\in I$ and $F$ with $\varphi_i(F)$ for some $i\in I$ (in fact, we have $\varphi_i(F)=\varphi_j(F)$ for all $i,j\in I$ as can easily be shown), we can thus assume without loss of generality $F\subseteq E$ and even $E_i\subseteq E$ for all $i\in I$ as all studied objects can also be studied over isomorphic fields. Thus, it is justified to speak about field extensions when it comes to direct limits. We will often use the above implicitly to simplify the notation.
\end{Bemerkung}

As a first step we will provide several examples of admissible field extensions that will be used later.

\begin{Lemma}\label{AnisotropicAdmissible}
	Let $X$ be a set of anisotropic quadratic forms over $F$ and $\mathcal C$ be the class of field extensions of $F$ such that every form $\varphi\in X$ stays anisotropic. Then $\mathcal C$ is admissible.
\end{Lemma}
\begin{Beweis}
	As $X$ consists of anisotropic quadratic forms we have $F\in\mathcal C$, thus (AD1). To verify (AD2) let $E=\varinjlim E_i$ with $E_i\in\mathcal C$ and $\varphi\in X$. If $\varphi_E$ was isotropic with isotropic vector $x\in E^{\dim\varphi}$, this vector $x$ would already be defined over some $E_j$ with $j\in I$. Thus $\varphi$ would be isotropic over $E_j$, a contradiction to $E_j\in\mathcal C$.\\
	Lastly if $E$ is a field extension of $F$ such that $\varphi_E$ is isotropic for some $\varphi\in X$, then $\varphi_K$ is clearly isotropic for every field extension $K/E$. Since this is equivalent to (AD3), the conclusion follows.
\end{Beweis}

\begin{Korollar}\label{RealAdmissible}
	The class of all formally real field extensions of $F$ is admissible.
\end{Korollar}
\begin{Beweis}
	This follows by applying \ref{AnisotropicAdmissible} to the set $X=\{n\times\qf1\mid n\in\N\}.$
\end{Beweis}

\begin{Lemma}\label{OrderingsAdmissible}
	Let $F$ be a real field and $X\subseteq X_F$ be a subset of orderings of $F$. Then the class $\mathcal C$ of field extensions extending the orderings in $X$ is admissible. In particular the class of totally positive field extensions of $F$ is admissible.
\end{Lemma}
\begin{Beweis}
	We clearly have $F\in\mathcal C$ so we have (AD1). Further \ref{ExtendOrderingsDirectLimit} implies the validity of (AD2).
	(AD3) holds due to the fact that the restriction of an ordering is again an ordering. \\
	To show that the class of totally positive field extensions of $F$ is admissible, we just have to apply the above to the case $X=X_F$.
\end{Beweis}

\begin{Bemerkung}
	Of course \ref{OrderingsAdmissible} can also be used to imply \ref{RealAdmissible}.
\end{Bemerkung}

\begin{Lemma}
	Let $F$ be a real field and $\mathcal C$ be the class of hereditarily totally positive field extensions of $F$. Then $\mathcal C$ is admissible.
\end{Lemma}
\begin{Beweis}
	We clearly have $F\in\mathcal C$, thus (AD1).
	Further (AD3) holds in the light of \ref{TotallyPositiveFieldTower} \ref{TotallyPositiveFieldTowerB}.\\
	For (AD2) let $(E_i)_{i\in I}$ be a directed system of fields in $\mathcal C$ and $E$ its direct limit in the category of field extensions of $F$. Let $K$ be an intermediate field, i.e. we have a tower of field extensions $E/K/F$. Let $P\in X_K$ be an ordering of $K$ and $a_1,\ldots, a_n\in P$ for some $n\in\N$. Let $x_1,\ldots, x_n\in E$ be such that $a_1x_1^2+\ldots+ a_nx_n^2=0$. Then, there is some $i\in I$ such that we have $a_k, x_k\in E_i$ for all $k\in\{1,\ldots, n\}$.\\
	We now consider the field $K_i:=E_i\cap K$. This field is ordered by $K_i\cap P$. As $E_i/F$ is an hereditarily totally positive extension by assumption, we can extend this ordering of $K_i$ to an ordering $P_i$ of $E_i$. Thus, the form $\qf{a_1,\ldots, a_n}$ is anisotropic over $E_i$ by \cite[VIII. Corollary 9.8]{Lam2005}, i.e. we have $x_1=\ldots=x_n=0$. Thus, the form is even anisotropic over $E$. Using \cite[VIII. Corollary 9.8]{Lam2005} again, we see that $P$ has an extension to $E$ and hence $E/K$ is totally positive. As $K$ was an arbitrary subfield of $E$, the extension $E/F$ is hereditarily totally positive and the proof is complete.
\end{Beweis}

The next result shows that we can combine the requirements of different admissible classes of field extensions to get a new admissible class.

\begin{Lemma}\label{SchnittAdmissible}
	Let $I$ be a nonempty index set and $\mathcal C_i$ be an admissible class of field extensions of the field $F$ for any $i\in I$. Then $\displaystyle\bigcap_{i\in I}\mathcal C_i$ is admissible.
\end{Lemma}
\begin{Beweis}
	As (AD1) and (AD3) together imply $F\in \mathcal C_i$ for all $i\in I$ the class $\displaystyle\bigcap_{i\in I}\mathcal C_i$ is not empty. The validity of (AD2) and (AD3) is clear.
\end{Beweis}

The crucial step to construct fields with supreme torsion forms is the following: 

\begin{Theorem}{\cite[6.1 Theorem]{Becher2004}}\label{BecherAdmissible}
	Let $\mathcal C$ be an admissible class of field extensions of $F$. There exists a field $K\in\mathcal C$ such that $K(\varphi)\notin\mathcal C$ for any anisotropic quadratic form $\varphi$ over $K$ of dimension at least 2.
\end{Theorem}

The proof of the above result uses a modification of Merkurjev's function field techniques. To be more precise, the field $K$ in \ref{BecherAdmissible} is constructed as a direct limit of iterated function field extensions. This fact will be used in an upcoming incidental remark. \\

A well known variation of the $u$-invariant of a field $F$ that we need to recall before stating the next result is the \textit{Hasse number} defined as 
$$\tilde u(F):=\max\{\dim\varphi\mid\varphi\text{ is a totally indefinite form over }F\}$$
or $\infty$ if no such maximum exists.

\begin{Proposition}\label{PropConstruction}
	Let $n\in\N$ be an integer with $n\geq2$ and $P$ be a non-empty set of anisotropic torsion $n$-fold Pfister forms over $F$. There exists a field extension $K/F$ such that
	\begin{enumerate}[(a)]
		\item\label{PropConstructionA} the field $K$ is formally real.
		\item\label{PropConstructionB} for every $\pi\in P$ the form $\pi_K$ is anisotropic .
		\item\label{PropConstructionC} we have $\tilde u(K)=u(K)=2^n$. In particular $I_t^{n+1}K$ is trivial.
		\item\label{PropConstructionD} every anisotropic form over $K$ that is indefinite at at least one ordering of $K$ is a subform of $\pi_K$ for some $\pi\in P$. 
		\item any anisotropic $n$-fold Pfister form over $K$ is 2-real maximal over K.
	\end{enumerate}
\end{Proposition}
\begin{Beweis}
	Let $\mathcal C$ be the class of field extensions $K/F$ such that $K$ is real and $\pi_K$ is anisotropic for every $\pi\in P$. This class is admissable due to \ref{AnisotropicAdmissible}, \ref{RealAdmissible} and \ref{SchnittAdmissible}.
	By \ref{BecherAdmissible} we get the existence of a field $K\in\mathcal C$ such that $K(\varphi)\notin\mathcal C$ for every anisotropic quadratic form $\varphi$ over $K$ of dimension at least 2.\\
	As $K\in\mathcal C$, it satisfies \ref{PropConstructionA} and \ref{PropConstructionB}.\\
	Now let $\varphi$ be a quadratic form that is indefinite at at least one ordering of $K$. Then $K(\varphi)$ is a real field by \ref{BspTotallyPositive} (c). Thus $\pi_{K(\varphi)}$ has to be isotropic and as a Pfister form therefore hyperbolic for some $\pi\in P$. By \cite[Chapter X. Corollary 4.9]{Lam2005} there exists some $a\in F^\ast$ with $a\varphi\subseteq\pi$. In particular $\dim\varphi\leq\dim\pi$, so we get \ref{PropConstructionC}. This further implies that $\pi_K$ is universal for every $\pi\in K$. Since a form is similar to a subform of a universal round form iff it is a subform of this form itself, we get \ref{PropConstructionD}.\\
	Now let $\pi$ be any anisotropic torsion $n$-fold Pfister form over $K$ and $a\in\sum F^{\ast2}\setminus F^{\ast2}$. As $n\geq 2$, the torsion form $\psi:=\qf{a}\perp\pi'$ is not similar to any Pfisterform and therefore has to be isotropic. That means $-a\in D_K(\pi')$, so we have $\qf{1,-a}\subseteq\pi$. The assertion follows by \ref{Connemara}.
\end{Beweis}

The next result can be proved in a similar way as \ref{PropConstruction} using a slightly modified admissible class of field extensions. We therefore will not give a complete proof but only the class of field extensions that can be used. We leave the details to the reader.

\begin{Proposition}\label{PropConstructionAlternative}
	Let $n\in\N$ be an integer with $n\geq2$ and $P$ be a set of anisotropic torsion $n$-fold Pfister forms over $F$. There exists a field extension $K/F$ such that
	\begin{enumerate}[(a)]
		\item The field extension $K/F$ is totally positive.
		\item For every $\pi\in P$ the form $\pi_K$ is anisotropic .
		\item We have $\tilde u(K)=u(K)=2^n$. In particular $I_t^{n+1}K$ is trivial.
		\item Every anisotropic totally indefinite form over $K$ is a subform of $\pi_K$ for some $\pi\in P$. 
		\item Any anisotropic $n$-fold Pfister form over $K$ is 2-real maximal over $K$.
	\end{enumerate}
\end{Proposition}
\begin{Beweis}
	This follows by considering the class 
	$$\mathcal C:=\{K/F\mid K/F\text{ is totally positive}, \pi_k\text{ is anisotropic for every }\pi\in P \},$$
	which is admissible by \ref{AnisotropicAdmissible}, \ref{OrderingsAdmissible} and \ref{SchnittAdmissible}.
\end{Beweis}

\begin{Bemerkung}
	We can give another proof for \ref{Construct2RealMaximal} using the just developed techniques. To do so, given an anisotropic torsion Pfister form $\varphi$ over the field $F$, consider the class of field extensions 
	$$\mathcal C:=\{E/F\mid E/F \text{ is a hereditarily totally positive extension, }\varphi \text{ anisotropic} \}.$$
	As the proof of \ref{BecherAdmissible} just uses the direct limit of several function field extensions, with \ref{HereditarilyTotallyRealFunctionField} in mind, one readily sees that the resulting field will be a 2-extension.
\end{Bemerkung}

\begin{Korollar}\label{ConstSTF}
	Let $\pi\in I_t^nF$ for some $n\in\N$ with $n\geq2$. Then there is a formally real field extension $K/F$ such that $\pi_K$ is the supreme torsion form over $K$. 
\end{Korollar}
\begin{Beweis}
	This is a direct consequence of \ref{PropConstruction} or \ref{PropConstructionAlternative} with ${P=\{\pi\}}$. 
\end{Beweis}

With the next example we will finish the discussion at the end of Section 3.

\begin{Beispiel}\label{2RealMaximalNotTorsion}
	We consider the form $\varphi=\qf{1,1, -3}$ over the field $\Q$ of rational numbers. It is clear that $\varphi$ is anisotropic and totally indefinite, but not a torsion form. We further consider the class of field extensions
	$$\mathcal C:=\{K/\Q\mid K/\Q\text{ is totally positive, }\varphi_K\text{ is anisotropic}\}.$$
	This class is admissible by \ref{AnisotropicAdmissible}, \ref{OrderingsAdmissible} and \ref{SchnittAdmissible}. By \ref{BecherAdmissible} there is some field $F\in\mathcal C$ such that $F(\psi)\notin\mathcal C$ for every anisotropic quadratic form $\psi$ over $F$ of dimension at least 2. Note that $F$ cannot be pythagorean as  otherwise $\varphi_F$ would be isotropic. Thus there is some $a\in\sum F^{\ast2}\setminus F^{\ast2}$. We then have $F(\pi)\notin\mathcal C$, where $\pi$ is the binary torsion form $\pi=\qf{1,-a}$. Using \ref{TotallyPositiveFieldTower} (\ref{TotallyPositiveFieldTowerA}), \ref{LemTotallyPositiveQuadExt} and the isomorphism $F(\pi)\cong F(\sqrt a)$, we see that $F(\pi)/\Q$ is totally positive. But as we have $F(\pi)\notin\mathcal C$, the form $\varphi_{F(\pi)}$ has to be isotropic which means that $\pi$ is similar to a subform of $\varphi_F$. In view of \ref{Connemara}, we see that we have constructed a field $F$ such that $\varphi$ is 2-real maximal over $F$, but $\varphi_F$ is not a torsion form.
\end{Beispiel}

\end{section}

\begin{section}{Correlations with Invariants}

While studying fields it is a natural task to find correlations between invariants of the given field. Several results in \cite{Pfister} deal with this question. One main problem is to find lower bounds for the number of square classes $|F^\ast/F^{\ast2}|$ in dependence of other invariants. In \cite{Pfister} this is done using the level $s(F)$ for nonreal fields and, for formally real fields, using the Pythagoras number. For a given formally real field, in order to show the estimate
$$|F^\ast/F^{\ast2}|\geq 2\cdot2^{\frac{t(t+1)}2}$$
for $t:=\lfloor\log_2p(F)\rfloor$, the idea in \cite[Satz 25]{Pfister} was to find lower bounds for the quotients $[D_F(2^n):D_F(2^{n-1})]$ for $n\in\{1,\ldots, t+1\}$ where $D_F(k):=\{x\in F^\ast\mid x\text{ is a sum of }k\text{ squares}\}$.
Pfister showed that we have $[D_F(2^{t+1}):D_F(2^t)]\geq2$. The main purpose in this section is to give an example to show that this bound is sharp. In order to have a more compact and consistent notation during the next results, we set $D_F(\infty)=\cup_{n\in\N}D_F(n)=\sum F^{\ast2}$.

\begin{Proposition}\label{PropReprElements}
	Let $n\in\N$ be an integer $\varphi\in P_nF$ the supreme torsion form over $F$ and $p(F)>2^{n-1}$. Then we have $[D_F(\infty):D_F(2^{n-1})]=2$.
\end{Proposition}
\begin{Beweis}
	As we have $p(F)>2^{n-1}$, we clearly have $[D_F(\infty):D_F(2^{n-1})]\geq2$, see also \cite[Satz 25, proof of Satz 18 d)]{Pfister}. So let now $x,y$ be representives of nontrivial classes of $D_F(\infty)/D_F(2^{n-1})$. Then ${2^{n-1}\times\qf1\perp\qf{-x}}$ and $2^{n-1}\times\qf1\perp\qf{-y}$ are anisotropic Pfister neighbors of the $n$-fold Pfister forms $2^{n-1}\times\Pfister x, 2^{n-1}\times\Pfister y$ respectively which therefore both have to be isometric to $\varphi$. Witt cancellation yields $2^{n-1}\times\qf{ -x}\cong2^{n-1}\times\qf{-y}$ which is equivalent to 
	$$2^{n-1}\times\qf{1}\cong2^{n-1}\times\qf{xy}.$$
	Thus we obtain $xy\in D_F(2^{n-1})$ which means that $x$ and $y$ represent the same class in $D_F(\infty)/D_F(2^{n-1})$.
\end{Beweis}

\begin{Bemerkung}
	The above proof shows in particular that in the situation of \ref{PropReprElements} the supreme torsion form is given by
	$$\Pfister{x,-1,\ldots, -1}$$
	for an arbitrary element $x\in D_F(\infty)\setminus D_F(2^{n-1})$, i.e. an element of length greater than $2^{n-1}$. Such forms will reoccur in \ref{ThmHeight}
\end{Bemerkung}

\begin{Beispiel}\label{BeispielPythagorasSupreme}
	To construct a field $K$ fulfilling the assumptions of \ref{PropReprElements} for a given $n\in\N$ with $n\geq2$ we can start with any field with pythagoras number greater than $2^{n-1}$. Specifically, we can choose $\R(X_1,\ldots, X_{2^{n-1}})$ by \cite[IX. Corollary 2.4]{Lam2005}. We then take an element $x$ of finite length greater than $2^{n-1}$. In the concrete example we can take $x=1+X_1^2+\ldots+X_{2^{n-1}}^2$. Finally we get the desired field $K$ by applying \ref{ConstSTF} to the given field and the $n$-fold Pfister form $\pi:=\Pfister{x, -1,\ldots, -1}$ which by construction is obviously an anisotropic torsion form.\\
	As $\pi_K$ is anisotropic, the element $x$ has to be of length greater than $2^{n-1}$ which implies $p(K)>2^{n-1}$, as desired.
\end{Beispiel}

\begin{Korollar}
	Let $F$ be a real field with finite Pythagoras number $p(F)$ fulfilling $2^{n-1}<p(F)\leq 2^{n}$ for some $n\in\N$. We then have $[D_F(\infty):D_F(2^{n-1})]\geq2$ and there are fields for which equality holds.
\end{Korollar}
\begin{Beweis}
	The estimate is clear and due to Pfister as mentioned above and fields for which the estimate is an equality are constructed in \ref{BeispielPythagorasSupreme}.
\end{Beweis}

As another invariant, we consider the \textit{height} of $F$, in symbols $h(F)$, i.e. the exponent of $W_tF$. As we only deal with formally real fields, this is the smallest 2-power $2^k$ with $2^k\geq p(F)$ if $p(F)\in\N$ or infinity otherwise, see \cite[Chapter XI. Theorem 5.6 (1)]{Lam2005}. The existence of supreme torsion forms has influence on the height as we see in the final result of this chapter.

\begin{Theorem}\label{ThmHeight}
	Let $F$ be a field with a supreme torsion form $\pi\in P_n(F)$ for some $n\in\N$. We then have $h(F)\leq 2^n$. We have equality if and only if $2^{n-1}\times\qf1\subseteq\pi$.
\end{Theorem}
\begin{Beweis}
	For the upper bound it is enough to show $p(F)\leq 2^n$. So let now $x$ be a non-zero sum of squares. The torsion Pfister form $2^n\times\qf{1,-x}$ cannot be a subform of the supreme torsion form due to the dimensions of the respective forms and thus has to be isotropic, hence even hyperbolic. Thus its Pfister neighbor $2^n\times\qf1\perp\qf{-x}$ is isotropic, which means that $x$ is a sum of at most $2^n$ squares as desired.\\
	For the equivalence, we start with the case $h(F)=2^n$. We then have $p(F)>2^{n-1}$, i.e. there is some $x\in \sum F^{\ast2}$ that is not a sum of $2^{n-1}$ squares. Thus the Pfister neighbor $2^{n-1}\times\qf1\perp\qf{-x}$ is anisotropic and therefore so is its associated Pfister form $2^{n-1}\times\qf{1,-x}$. The latter form therefore is an anisotropic $n$-fold torsion Pfister form and has to be the supreme torsion form $\pi$. As this forms has $2^{n-1}\times\qf{1}$ as a subform, this implication is done.\\
	For the other implication, we now assume $2^{n-1}\times\qf1\subseteq\pi$. Let $x\in F^\ast$ such that $2^{n-1}\times\qf{1}\perp\qf{-x}\subseteq\pi$. By considering the signature of this form and using Artins Theorem, we obtain ${x\in\sum F^{\ast2}}$. But $x$ cannot be a sum of at most $2^{n-1}$ squares because otherwise, the above form and thus also $\pi$ would be isotropic. This implies $2^{n-1}<p(F)$. With $h(F)\leq 2^n$ and the description of the height in terms of the Pythagoras number given above, we are done.
\end{Beweis}

\end{section}

\section*{Acknowledgments and Notes}
The results contained in this paper are part of the PhD-Thesis of the author. He would like to thank Detlev Hoffmann for supervising this work and giving some very useful hints and corrections in the process.

\bibliographystyle{abbrv}
\bibliography{literatur}

\begin{thebibliography}{1}

\bibitem{Becher2004}
K.~J. Becher.
\newblock Supreme {P}fister forms.
\newblock {\em Comm. Algebra}, 32(1):217--241, 2004.

\bibitem{Becher2006}
K.~J. Becher.
\newblock Totally positive extensions and weakly isotropic forms.
\newblock {\em Manuscripta Math.}, 120(1):83--90, 2006.

\bibitem{ElmanKarpenkoMerkurjev2008}
R.~Elman, N.~Karpenko, and A.~Merkurjev.
\newblock {\em The algebraic and geometric theory of quadratic forms},
  volume~56 of {\em American Mathematical Society Colloquium Publications}.
\newblock American Mathematical Society, Providence, RI, 2008.

\bibitem{Lam2005}
T.~Y. Lam.
\newblock {\em Introduction to quadratic forms over fields}, volume~67 of {\em
  Graduate Studies in Mathematics}.
\newblock American Mathematical Society, Providence, RI, 2005.

\bibitem{Morandi}
P.~Morandi.
\newblock {\em Field and {G}alois theory}, volume 167 of {\em Graduate Texts in
  Mathematics}.
\newblock Springer-Verlag, New York, 1996.

\bibitem{Pfister}
A.~Pfister.
\newblock Quadratische {F}ormen in beliebigen {K}\"{o}rpern.
\newblock {\em Invent. Math.}, 1:116--132, 1966.

\bibitem{Prestel}
A.~Prestel.
\newblock {\em Lectures on formally real fields}, volume 1093 of {\em Lecture
  Notes in Mathematics}.
\newblock Springer-Verlag, Berlin, 1984.

\bibitem{ZalesskiiRibes2000}
L.~Ribes and P.~Zalesskii.
\newblock {\em Profinite groups}, volume~40 of {\em Ergebnisse der Mathematik
  und ihrer Grenzgebiete. 3. Folge. A Series of Modern Surveys in Mathematics
  [Results in Mathematics and Related Areas. 3rd Series. A Series of Modern
  Surveys in Mathematics]}.
\newblock Springer-Verlag, Berlin, 2000.

\bibitem{Scharlau}
W.~Scharlau.
\newblock {\em Quadratic and {H}ermitian forms}, volume 270 of {\em Grundlehren
  der Mathematischen Wissenschaften [Fundamental Principles of Mathematical
  Sciences]}.
\newblock Springer-Verlag, Berlin, 1985.

\end{thebibliography}

\end{document}